\documentclass[MSNbibl,number,citesort,dvips]{arxbj}
\usepackage{upgreek}
\usepackage{graphicx}

\aid{0}
\volume{18}
\issue{4}
\pubyear{2012}
\firstpage{1405}
\lastpage{1420}
\doi{10.3150/11-BEJ381}

\makeatletter
\newtheorem{theo}{Theorem}
\newremark{remark}{Remark}
\newcommand{\bbA}{{\mathbf A}}
\newcommand{\bbB}{{\mathbf B}}
\newcommand{\bbI}{{\mathbf I}}
\newcommand{\bbS}{{\mathbf S}}
\newcommand{\bbX}{{\mathbf X}}
\newcommand{\bbz}{{\mathbf z}}
\newcommand{\binoma}[2]{{#1\choose #2}}
\newcommand{\binom}[2]{\pmatrix{#1\cr #2}}
\newcommand{\fraca}[2]{{#1}/{#2}}
\makeatother

\begin{document}
\begin{frontmatter}

\title{Convergence of the largest eigenvalue of normalized sample
covariance matrices when $p$ and $n$ both tend to infinity with
their ratio converging to zero}
\runtitle{Convergence of the largest eigenvalue of normalized sample
covariance matrices}

\begin{aug}
\author{\fnms{B.B.} \snm{Chen}\corref{}\thanksref{e1}\ead[label=e1,mark]{chen0635@e.ntu.edu.sg}} \and
\author{\fnms{G.M.} \snm{Pan}\thanksref{e2}\ead[label=e2,mark]{gmpan@ntu.edu.sg}}
\runauthor{B.B. Chen and G.M. Pan} 
\address{Division of Mathematical Sciences, School of Physical and
Mathematical Science, Nanyang Technological University, Singapore.
\printead{e1,e2}}
\end{aug}

\received{\smonth{5} \syear{2010}}
\revised{\smonth{12} \syear{2010}}

%
\begin{abstract}
Let $\bbX_{p}=({\mathbf s}_1,\ldots,{\mathbf s}_n)=(X_{ij})_{p\times n}$ where
$X_{ij}$'s are independent and identically distributed (i.i.d.)
random variables with $EX_{11}=0,EX_{11}^{2}=1$ and
$EX_{11}^{4}<\infty$. It is showed that the largest eigenvalue of
the random matrix
$\bbA_{p}=\frac{1}{2\sqrt{np}}(\bbX_{p}\bbX_{p}^{\prime}-n\bbI_{p})$
tends to $1$ almost surely as
$p\rightarrow\infty,n\rightarrow\infty$ with $p/n \rightarrow0$.
\end{abstract}

%
\begin{keyword}
\kwd{empirical distribution}
\kwd{maximum eigenvalue}
\kwd{random matrices}
\end{keyword}

\end{frontmatter}

\section{Introduction}

Consider the sample covariance type matrix
$\bbS=\frac{1}{n}\bbX_{p}\bbX_{p}^{\prime}$, where
$\bbX_{p}=({\mathbf s}_1,\ldots,{\mathbf s}_n)=(X_{ij})_{p\times n}$ and
$X_{ij},i=1,\ldots,p,j=1,\ldots,n$, are i.i.d. random variables with
mean zero and variance 1. For such a matrix, much attention has been
paid to asymptotic properties of its eigenvalues in the setting of
$p/n\rightarrow c>0$ as $p\rightarrow\infty$ and
$n\rightarrow\infty$. For example, its empirical spectral
distribution (ESD) function $F^{\bbS}(x)$ converges with probability
one to the famous Mar\v{c}enko and Pastur law (see~\cite{MP} and
\cite{Jon}). Here, the ESD for any matrix $\bbA$ with real
eigenvalues $\lambda_{1}\leq\lambda_{2}\leq\cdots \leq\lambda_{p} $ is
defined by
\[
F^{\bbA}(x)=\frac{1}{p}\#\{i\dvt\lambda_{i}\leq x\},
\]
where $\#\{\cdots \}$ denotes the number of elements of the set. Also,
with probability one its maximum eigenvalue and minimum eigenvalue
converge, respectively, to the left end point and right end point of
the support of Mar\v{c}enko and Pastur's law (see~\cite{Ger} and
\cite{b3}).

In contrast with asymptotic behaviors of $\bbS$ in the case of
$p/n\rightarrow c$, the asymptotic properties of $\bbS$ have not
been well understood when $p/n\rightarrow0$. The first breakthrough
was made in Bai and Yin~\cite{b1}. They considered the normalized matrix
\[
\bbA_{p}=\frac{1}{2\sqrt{np}}(\bbX_{p}\bbX_{p}^{\prime}-nI_{p})
\]
and proved with probability one
\[
F^{\bbA_{p}}\rightarrow F(x),
\]
which is the so-called semicircle law with a density
\begin{eqnarray*}
\label{6} F'(x)=\left\{
\everymath{\displaystyle }
\begin{array}{l@{ \qquad }l}
\frac{2}{\uppi}\sqrt{1-x^{2}}, & \mbox{if } |x| \leq1,\\[5pt]
0, & \mbox{if } |x| > 1.
\end{array}
\right. \label{6}
\end{eqnarray*}
One should note that the semicircle law is also the limit of the
empirical spectral distribution of a symmetric random matrix whose
diagonal are i.i.d. random variables and above diagonal elements are
also i.i.d. (see~\cite{w}). Second, when $X_{11}\sim N(0,1)$, El
Karoui~\cite{ka0} proved that the largest eigenvalue of $\bbX_p\bbX_p'$
after properly centering and scaling converges to the Tracy$-$Widom law.

In this paper, for general $X_{11}$, we investigate the maximum
eigenvalue of $\bbA_{p}$ under the
setting of $p/n\rightarrow0$ as $p\rightarrow\infty$ and
$n\rightarrow\infty$. The main results are presented
in the following theorems.

%
\begin{theo}\label{theo1}
Let $\bbX_p=(X_{ij})_{p\times n}$ where
$\{X_{ij}\dvt i=1,2,\ldots,p;j=1,2,\ldots,n\}$ are i.i.d. real random
variables with $EX_{11}=0,EX_{11}^{2}=1 $ and $EX_{11}^{4}<\infty$.
Suppose that $n=n(p)\to\infty$ and $p/n \rightarrow0$
as $p\rightarrow\infty$. Define
\[
\bbA_{p}=(A_{ij})_{p\times p}=\frac{1}{2\sqrt{np}}(\bbX_{p}\bbX
_{p}^{\prime}-nI_{p}).
\]
Then as $p\rightarrow\infty$
\[
\lambda_{\mathrm{max}}(\bbA_{p})\to1 \qquad    \mbox{a.s.},
\]
where $\lambda_{\mathrm{max}}(\bbA_{p})$ represents the largest eigenvalue of
$\bbA_{p}$.
\end{theo}

Indeed, after truncation and normalization of the entries of the
matrix $\bbA_p$, we may obtain a better result.

%
\begin{theo}\label{theo2}
Let $n=n(p)\rightarrow\infty$ and $p/n\rightarrow0$ as $p\rightarrow
\infty$. Define a $p\times p$ random matrix~$\bbA_p$:
\[
\bbA_{p}=(A_{ij})_{p\times p}=\frac{1}{2\sqrt{np}}(\bbX_{p}\bbX
_{p}^{\prime}-nI_{p}),
\]
where $\bbX_p=(X_{ij})_{p\times n}$. Suppose that $X_{ij}$'s are
i.i.d. real random variables and satisfy the following conditions
\begin{enumerate}[(2)]
\item[(1)] $EX_{11}=0,EX_{11}^{2}= 1,EX_{11}^{4}<\infty$ and
\item[(2)] $|X_{ij}|\leq\delta_{p}\sqrt[4]{np}$, where
$\delta_{p}\downarrow0$, but $\delta_{p}\sqrt[4]{np}\uparrow
+\infty$, as $p\rightarrow\infty$.
\end{enumerate}

Then, for any $\epsilon>0,\ell>0$
\[
p\bigl(\lambda_{\mathrm{max}}(\bbA_{p})\geq1+\epsilon\bigr)=\mathrm{o}(p^{-\ell}).
\]
\end{theo}

So far we have considered the sample covariance type matrix $\bbS$.
However, a common used sample covariance matrix in statistics is
\[
\bbS_1=\frac{1}{n}\sum_{j=1}^n({\mathbf s}_j-\bar{\mathbf s})({\mathbf
s}_j-\bar{\mathbf s})',
\]
where
\[
\bar{\mathbf s}=\frac{1}{n}\sum_{j=1}^n{\mathbf s}_j.
\]
Similarly
we renormalize it as
\[
\bbA_{p1}=\frac{1}{2}\sqrt{\frac{n}{p}}(\bbS_1-I_{p}).
\]

\begin{theo}\label{theo3}
Under assumptions of Theorem~\ref{theo1}, as $p\rightarrow\infty$
\[
\lambda_{\mathrm{max}}(\bbA_{p1})\to1   \qquad  \mbox{a.s.},
\]
where $\lambda_{\mathrm{max}}(\bbA_{p1})$ stands for the largest eigenvalues
of $\bbA_{p1}$.
\end{theo}

Estimating a population covariance matrix for high dimension data is
a challenging task. Usually, one can not expect the sample
covariance matrix to be a consistent estimate of a population
covariance matrix when both $p$ and $n$ go to infinity, especially
when the orders of $p$ and $n$ are very close to each other. In such
circumstance, as argued in~\cite{Ka}, operator norm consistent
estimation of large population covariance matrix still has nice
properties.

Suppose that $\Sigma$ is a population covariance matrix, nonnegative
definite symmetric matrix. Then $\Sigma^{1/2}{\mathbf
s}_j,j=1,\ldots,n$,
may be viewed as i.i.d. sample drawn from the population with
covariance matrix $\Sigma$, where $(\Sigma^{1/2})^2=\Sigma$. The
corresponding sample covariance matrix is
\[
\bbS_{2}=\frac{1}{n}\sum_{j=1}^n(\Sigma^{1/2}{\mathbf
s}_j-\Sigma^{1/2}\bar{\mathbf s})(\Sigma^{1/2}{\mathbf s}_j-\Sigma
^{1/2}\bar{\mathbf s})'.
\]
Theorem~\ref{theo3} indicates that the matrix $\bbS_{2}$ is an
operator consistent estimation of $\Sigma$ as long as
$p/n\rightarrow0$ when $p\rightarrow\infty$. Specifically, we have
the following theorem.

\begin{theo}\label{theo4}
In addition to the assumptions of Theorem~\ref{theo1}, assume that
$\|\Sigma\|$ is bounded. Then, as $p\rightarrow\infty$
\[
\|\bbS_{2}-\Sigma\|=\mathrm{O}\Biggl(\sqrt{\frac{p}{n}}\Biggr) \qquad    \mbox{a.s.},
\]
where $\|\cdot\|$ stands for the  spectral norm of a matrix.
\end{theo}

\begin{remark}
Related work is~\cite{raa}, where the authors investigated
quantitative estimates of the convergence of the empirical
covariance matrix in the Log-concave ensemble. Here we obtain a
convergence rate of the empirical covariance matrix when the sample
vectors are in the form of $\Sigma^{1/2}{\mathbf s}_j$.
\end{remark}

\begin{remark}
Theorems~\ref{theo1}--\ref{theo4} are stated for the real random matrix $\bbX_p$, but they
also hold for the complex case under moment conditions
$EX_{11}=0,E|X_{11}|^2=1$ and $E|X_{11}|^4<\infty$. The proofs are
similar to those for the real case except some notation changes.
\end{remark}

\section{\texorpdfstring{Proof of Theorem \protect\ref{theo1}}{Proof of Theorem 1}}
Throughout the paper, $C$ denotes a constant whose value may vary from
line to line. Also, all limits in the paper are taken as $p\to\infty
$.

It follows from Theorem in~\cite{b1} that
%
\begin{equation}\label{a1}
\liminf_{p\rightarrow\infty}\lambda_{\max}(\bbA_{p})\geq
1   \qquad \mbox{a.s.}
\end{equation}
Thus, it suffices to show that
%
\begin{equation}\label{a2}
\limsup_{p\rightarrow\infty}\lambda_{\max}(\bbA_{p})\leq
1  \qquad  \mbox{a.s.}
\end{equation}

Let $\hat{\bbA}_{p}=\frac{1}{2\sqrt{np}}(\hat{\bbX}_{p}\hat{\bbX
}_{p}^{\prime}-nI_{p})$, where $\hat{\bbX}_{p}=(\hat{X}_{ij})_{p\times n}$
and $\hat{X}_{ij}=X_{ij}I(|X_{ij}|\leq\delta_p\sqrt[4]{np})$ where
$\delta_p$ is chosen as the larger of $\delta_p$ constructed as in
(\ref{a3}) and $\delta_p$ as in (\ref{b1}). On the one hand, since
$EX_{11}^4<\infty$ for any $\delta>0$
we have
\[
\mathop{\lim}\limits_{p\to\infty}\delta
^{-4}E|X_{11}|^4I\bigl(|X_{11}|>\delta\sqrt[4]{np}\bigr)=0.
\]
Since the above is true for arbitrary positive $\delta$ there exists a
sequence of positive $\delta_p$ such that
%
\begin{equation}\label{a3}
\lim_{p\rightarrow\infty}\delta_p=0, \qquad
\mathop{\lim}\limits_{p\to\infty}\delta
_p^{-4}E|X_{11}|^4I\bigl(|X_{11}|>\delta_p\sqrt[4]{np}\bigr)=0, \qquad  \delta
_{p}\sqrt[4]{np}\uparrow
+\infty.
\end{equation}
On the other hand, since $EX_{11}^4<\infty$ for any $\nu>0$
\[
\sum_{k=1}^{\infty} 2^{k}
P (|X_{11}|>\nu2^{\fraca{k}{4}} )<\infty.
\]
In view of the arbitrariness of $\nu$, there is a sequence of
positive number $\nu_k$ such that
%
\begin{equation}\label{c1}
\nu_k\rightarrow0, \mbox{ as }  k\rightarrow\infty, \qquad
\sum_{k=1}^{\infty} 2^{k}
P (|X_{11}|>\nu_k2^{\fraca{k}{4}} )<\infty.
\end{equation}
For each $k$, let $p_k$ be the maximum $p$ such that
$n(p)\cdot p\leq2^k$. For $p_{k-1}<p\leq p_k$, set
%
\begin{equation}\label{b1}\delta_p=2\nu_k.
\end{equation}
Let $\label{c2}Z_t=X_{ij},t=(i-1)n+j$ and obviously $\{Z_t\}$ are
i.i.d.
We then conclude from (\ref{c1}) and (\ref{b1}) that
\begin{eqnarray*}
P(\bbA_p\neq\hat{\bbA}_{p},i.o.)&\leq&
\mathop{\lim}\limits_{K\to\infty}P \Biggl(\mathop{\bigcup}\limits
_{k=K}^{\infty}\mathop{\bigcup}\limits_{p_{k-1}<p\leq p_k}\mathop
{\bigcup}\limits_{i\leq p,j\leq n} \bigl\{|X_{ij}|>\delta_p\sqrt
[4]{np} \bigr\} \Biggr)\\[-2pt]
  &\leq&\mathop{\lim}\limits_{K\to\infty}\sum_{k=K}^{\infty}P
\Biggl(\mathop{\bigcup}\limits_{p_{k-1}<p\leq p_k}\mathop{\bigcup
}\limits_{t=1}^{2^{k}} \{|Z_t|>\nu_k2^k \}
\Biggr)\\[-2pt]
  &=&\mathop{\lim}\limits_{K\to\infty}\sum_{k=K}^{\infty}P
\Biggl(\mathop{\bigcup}\limits_{t=1}^{2^{k}} \{|Z_t|>\nu
_k2^k \} \Biggr)\\[-2pt]
  &\leq&\mathop{\lim}\limits_{K\to\infty}\sum_{k=K}^{\infty} 2^{k}
P (|Z_1|>\nu_k2^{\fraca{k}{4}} )\\[-2pt]
   &=& 0   \qquad  \mbox{a.s.}
\end{eqnarray*}
It follows that
$\lambda_{\max}(\bbA_p)-\lambda_{\max}(\hat{\bbA}_p)\to0$  a.s.
as $p\rightarrow\infty$.

From now on, we write $\delta$ for $\delta_p$ to simplify
notation.
Moreover, set
$\tilde{\bbA}_{p}=\frac{1}{2\sqrt{np}}(\tilde{\bbX}_{p}\tilde
{\bbX}_{p}^{\prime}-nI_{p})$,
where $\tilde{\bbX}_{p}=(\tilde{X}_{ij})_{p\times n}$ and $\tilde
{X}_{ij}=\frac{\hat{X}_{ij}-E\hat{X}_{11}}{\sigma}$. Here, $\sigma
^2=E(\hat{X}_{11}-E\hat{X}_{11})^2$ and $\sigma^2\to1$ as $p\to
\infty$.

We obtain via (\ref{a3})
%
\begin{equation}\label{a4}
|E\hat{X}_{11}|\leq\frac{E|X_{11}|^4I(|X_{11}|>
\delta_p\sqrt[4]{np})}{\delta^3(np)^{3/4}}\leq\frac{C}{(np)^{3/4}}
\end{equation}
and
%
\begin{equation}\label{a5}
|\sigma^2-1|\leq CE|X_{11}|^2I\bigl(|X_{11}|>\delta\sqrt[4]{np}\bigr)\leq
\frac{E|X_{11}|^4I(|X_{11}|>\delta\sqrt[4]{np})}{\delta^2\sqrt
{np}}=\mathrm{o}\biggl(\frac{1}{\sqrt{np}} \biggr).
\end{equation}

We conclude from the Rayleigh--Ritz theorem that
\begin{eqnarray*}
&&|\lambda_{\mathrm{max}}(\tilde{\bbA}_p)-\lambda_{\mathrm{max}}(\hat{\bbA}_p)|
\\[-2pt]
&& \qquad \leq\frac{1}{2\sqrt{np}} \Biggl|\mathop{\sup}\limits_{\|\bbz\|
=1}\Biggl (\sum_{i\neq j}z_iz_j\sum_{k=1}^{n}\hat{X}_{ik}\hat
{X}_{jk}+ \sum_{i=1}^{p}z_i^2\sum_{k=1}^{n}(\hat{X}_{ik}^2-1)
\Biggr)
\\[-2pt]
&&\hphantom{ \leq\frac{1}{2\sqrt{np}} \Biggl|} \qquad {}  -\mathop{\sup}\limits_{\|\bbz\|=1} \Biggl(\sum_{i\neq
j}z_iz_j\sum_{k=1}^{n}\tilde{X}_{ik}\tilde{X}_{jk}+\sum
_{i=1}^{p}z_i^2\sum_{k=1}^{n}(\tilde{X}_{ik}^2-1) \Biggr) \Biggr|
\\[-2pt]
&& \qquad  \leq\frac{1}{2\sqrt{np}}  \biggl|1-\frac{1}{\sigma^2}
\biggr|\mathop{\sup}\limits_{\|\bbz\|=1} \Biggl|\sum_{i\neq j}z_iz_j\frac
{1}{\sqrt{np}}\sum_{k=1}^{n}\hat{X}_{ik}\hat{X}_{jk}+ \sum
_{i=1}^{p}z_i^2\sum_{k=1}^{n}(\hat{X}_{ik}^2-1) \Biggr| \\[-2pt]
&&  \quad \qquad {}  +\frac{1}{2\sqrt{np}} \frac{2|EX_{11}|}{\sigma^2} \mathop
{\sup}\limits_{\|\bbz\|=1}  \Biggl| \sum_{i=1}^{p}\sum
_{j=1}^{p}z_iz_j\sum_{k=1}^{n}\hat{X}_{ik}  \Biggr|\\[-2pt]
&& \quad  \qquad {}  +\frac{1}{2\sqrt{np}} \frac{n|EX_{11}|^2}{\sigma^2} \mathop
{\sup}\limits_{\|\bbz\|=1}  \Biggl|\sum_{i=1}^{p}\sum
_{j=1}^{p}z_iz_j \Biggr| +\frac{n}{2\sqrt{np}}  \biggl|1-\frac
{1}{\sigma^2} \biggr|\\[-2pt]
&& \qquad  = A_1+A_2+A_3+A_4.
\end{eqnarray*}
By (\ref{a5}) and the strong law of large numbers, we have
\begin{eqnarray*}
A_1&=&\frac{|\sigma^2-1|}{2\sqrt{np}\sigma^2} \mathop{\sup}\limits
_{\|\bbz\|=1}  \Biggl| \sum_{k=1}^{n}\Biggl ( \Biggl(\sum
_{i=1}^{p}z_i\hat{X}_{ik}  \Biggr)^2-\sum_{i=1}^{p}z_i^2\hat
{X}_{ik}^2  \Biggr)+\sum_{i=1}^{p}z_i^2\sum_{k=1}^{n} (\hat
{X}_{ik}^{2}-1 )  \Biggr|\\
& \leq&\frac{|\sigma^2-1|\sqrt{np}}{2\sigma^2}\cdot\frac
{1}{np} \Biggl(2 \Biggl|\sum_{k=1}^{n}\sum_{i=1}^{p}\hat{X}_{ik}^2
 \Biggr|+  \Biggl|\sum_{i=1}^{p}\sum_{k=1}^{n} (\hat
{X}_{ik}^{2}-1 ) \Biggr|  \Biggr)\\
& \leq&\frac{|\sigma^2-1|\sqrt{np}}{2\sigma^2}\cdot\frac
{1}{np}\Biggl (3 \Biggl|\sum_{k=1}^{n}\sum_{i=1}^{p}X_{ik}^2
\Biggr|+np  \Biggr)\\
& \to&0   \qquad \mbox{a.s.}
\end{eqnarray*}
Similarly, (\ref{a4}), H\"{o}lder's inequality and the strong
law of large numbers yield
\begin{eqnarray*}
A_2& \leq&\frac{1}{2\sqrt{np}}\cdot\frac{2|E\hat{X}_{11}|}{\sigma
^2} \mathop{\sup}\limits_{\|\bbz\|=1} \Biggl|\sum
_{j=1}^{p}z_j \Biggr| \Biggl|\sum_{i=1}^{p}z_i\sum_{k=1}^{n}\hat
{X}_{ik} \Biggr|\\
& \leq&\frac{1}{2\sqrt{np}}\cdot\frac{C}{\sigma^2(np)^{3/4}}\cdot
\sqrt{p}\cdot \Biggl(\sum_{i=1}^{p} \Biggl(\sum_{k=1}^{n}\hat
{X}_{ik} \Biggr)^2 \Biggr)^{1/2}\\
&\leq&\frac{1}{2\sqrt{np}}\cdot\frac{C}{\sigma^2(np)^{3/4}}\cdot
\sqrt{p}\cdot\Biggl (n\sum_{i=1}^{p}\sum_{k=1}^{n}\hat
{X}_{ik}^2 \Biggr)^{1/2}
\\
&\leq&\frac{C}{\sigma^2(np)^{1/4}} \Biggl|\frac{1}{np}\sum
_{i=1}^{p}\sum_{k=1}^{n}\hat{X}_{ik}^2 \Biggr|^{1/2}
\\
&\leq&\frac{C}{\sigma^2(np)^{1/4}} \Biggl|\frac{1}{np}\sum
_{i=1}^{p}\sum_{k=1}^{n}X_{ik}^2 \Biggr|^{1/2}\to0   \qquad  \mbox{a.s.}
\end{eqnarray*}

It is straightforward to conclude from (\ref{a4}) and (\ref{a5}) that
\[
A_3\to0  \qquad  \mbox{a.s.},  \qquad  A_4\to0  \qquad  \mbox{a.s.}
\]

Thus, we have $\lambda_{\max}(\hat{\bbA}_p)-\lambda_{\max}(\tilde
{\bbA}_p)\to0$
a.s. By the above results, to prove (\ref{a2}), it is sufficient to
show that $\limsup_{p\rightarrow\infty}\lambda_{\max
}(\tilde{\bbA}_p)\leq1$  a.s. To this end, we note that the
matrix $\tilde{\bbA}_p$ satisfies all the assumptions in Theorem \ref
{theo2}. Therefore, we obtain (\ref{a2}) by Theorem~\ref{theo2}
(whose argument is given in the next section). Together with (\ref
{a1}), we finishes the proof of Theorem~\ref{theo1}.

\section{\texorpdfstring{Proof of Theorem \protect\ref{theo2}}{Proof of Theorem 2}}
Suppose that $\bbz=(z_1,\ldots,z_p)$ is a unit vector. By the
Rayleigh--Ritz theorem, we then have
\begin{eqnarray*}
\lambda_{\max}(\bbA_{p}) & = &\mathop{\max}\limits_{\parallel z
\parallel=1}\biggl(\sum_{i,j}z_{i}z_{j}A_{ij}\biggr) \\
& = &\mathop{\max}\limits_{\parallel z \parallel=1}\Biggl(\sum_{i\neq
j}z_{i}z_{j}A_{ij}+\sum_{i=1}^{p}z_{i}^{2}A_{ii}\Biggr)\\
&\leq&\lambda_{\max}(\bbB_{p})+\mathop{\max}\limits_{i\leq
p}|A_{ii}|,
\end{eqnarray*}
where $\bbB_p=(B_{ij})_{p\times p}$ with
\[
B_{ij}=\left\{
\everymath{\displaystyle }
\begin{array}{l@{ \qquad }l}
0, & \mbox{if } i=j,\\
\frac{1}{2\sqrt{np}}\sum_{k=1}^{n}X_{ik}X_{jk}, & \mbox{if } i\neq
j.
\end{array}
\right.
\]

To prove Theorem~\ref{theo2}, it is sufficient to prove, for any $\epsilon
>0,\ell>0$
%
\begin{equation}\label{a6}
P \bigl(\lambda_{\mathrm{max}}(\bbB_{p})>1+\epsilon \bigr)=\mathrm{o}(p^{-l})
\end{equation}
and
%
\begin{equation}\label{a7}
P\Biggl(\mathop{\max}\limits_{i\leq p}\frac{1}{\sqrt{np}}\Biggl|\mathop{\sum
}\limits_{j=1}^{n}(X_{ij}^2-1)\Biggr|>\epsilon\Biggr)=\mathrm{o}(p^{-l}).
\end{equation}

We first prove (\ref{a7}). To simplify notation, let $Y_j=X_{1j}^2-1 $
and $ C_1=E|Y_1|^2$. Then $EY_j=0$. Choose an appropriate sequence
$h=h_{p}$ such that it satisfies, as $p\to\infty$
%
\begin{equation}\label{a8}
\left\{
\everymath{\displaystyle }
\begin{array}{l}
  h/\log{p}\to\infty,\\
  \delta^{2}h/\log{p}\to0,\\
  \frac{\delta^{4}p}{C_1}\geq\sqrt{p}.
\end{array}
\right.
\end{equation}
We then have
\begin{eqnarray*}
&& P\Biggl(\mathop{\max}\limits_{i\leq p}\frac{1}{\sqrt{np}}\Biggl|\mathop
{\sum}\limits_{j=1}^{n}(X_{ij}^2-1)\Biggr|>\epsilon\Biggr)\\
&& \quad \leq p\cdot P
\Biggl(\Biggl|\sum_{j=1}^{n}(X_{1j}^2-1)\Biggr|>\epsilon\sqrt{np} \Biggr)\\
&& \quad  \leq\epsilon^{-h}p\bigl(\sqrt{np}\bigr)^{-h}E\Biggl|\sum_{j=1}^{n}Y_j\Biggr|^h
\\
&& \quad\leq\epsilon^{-h}p\bigl(\sqrt{np}\bigr)^{-h}\sum_{m=1}^{h/2}\sum_{1\leq
j_1<j_2<j_m\leq n}\mathop{\mathop{\sum}_{ i_1+i_2+\cdots +i_m=h }}_{i_1\geq2,\ldots,i_1\geq2 }\frac{h!}{i_1!i_2!\cdots i_m!}E|Y_{j_1}|^{i_1}E|Y_{j_2}|^{i_2}\cdots E|Y_{j_m}|^{i_m}\\
&& \quad\leq\epsilon^{-h}p\bigl(\sqrt{np}\bigr)^{-h}\sum_{m=1}^{h/2}\mathop{\mathop{\sum
}_{   i_1+i_2+\cdots +i_m=h }}_{i_1\geq2,\ldots,i_1\geq2   }\frac{n!}{m!(n-m)!}\frac
{h!}{i_1!i_2!\cdots i_m!}E|Y_1|^{i_1}E|Y_1|^{i_2}\cdots E|Y_1|^{i_m}\\
&& \quad \leq\epsilon^{-h}p\bigl(\sqrt{np}\bigr)^{-h}\sum_{m=1}^{h/2}\mathop{\mathop{\sum
}_{ i_1+i_2+\cdots +i_m=h }}_{ i_1\geq2,\ldots,i_1\geq2  } n^m\frac
{h!}{i_1!i_2!\cdots i_m!}C_1^m\bigl(\delta^2\sqrt{np}\bigr)^{h-2m}\\
&& \quad \leq\epsilon^{-h}p\sum_{m=1}^{h/2}m^h \biggl(\frac{\delta
^4p}{C_1} \biggr)^{-m}\delta^{2h} \leq\epsilon^{-h}p\frac
{h}{2}\cdot \biggl(\frac{\delta^2h}{\log{(\delta^4p/C_1)}}
\biggr)^h\\
&& \quad \leq \biggl( \biggl(\frac{ph}{2} \biggr)^{1/h}\cdot\frac{2\delta
^2h}{\log{p}}\cdot\epsilon^{-1} \biggr)^h \leq \biggl(\frac{\xi
}{\epsilon} \biggr)^h=\mathrm{o}(p^{-\ell}),
\end{eqnarray*}
where $\xi$ is a constant satisfying $0<\xi<\epsilon$. Below are
some interpretations of the above inequalities:
\begin{longlist}[(b)]
\item[(a)] The fifth inequality is because, $\frac{n!}{m!(n-m)!}<n^m$,
$|Y_1|<\delta^2\sqrt{np}$.
\item[(b)] We use the fact $\mathop{\mathop{\sum}_{
i_1+i_2+\cdots +i_m=h }}_{ i_1\geq2,\ldots,i_1\geq2  }\frac{h!}{i_1!i_2!\cdots i_m!}<m^h$ in the sixth inequality.
\item[(c)] The seventh inequality uses the elementary inequality
\[
a^{-t}t^b\leq \biggl(\frac{b}{\log{a}} \biggr)^b  \qquad \mbox{for
all }  a>1,b>0,t\geq1 \mbox{ and }  \frac{b}{\log{a}}>1.
\]
\item[(d)] The last two inequalities are due to (\ref{a8}).
\item[(e)] With the facts that $\frac{\xi}{\epsilon}<1,h/\log{p}\to
\infty$, the last equality is true.
\end{longlist}
Thus, (\ref{a7}) follows.

Next, consider (\ref{a6}). For any $\varsigma>0$, we have
\begin{eqnarray*}
 P\bigl(\lambda_{\mathrm{max}}(\bbB_{p})\geq1+\varsigma\bigr)&\leq&
\frac{E\lambda_{\mathrm{max}}^{k}(\bbB_{p})}{(1+\varsigma)^{k}} \leq\frac
{Etr(\bbB_{p}^{k})}{(1+\varsigma)^{k}}\\
& \leq&
\frac{1}{(1+\varsigma)^{k}}\cdot\frac{1}{(2\sqrt{np})^{k}}\sum
E(X_{i_{1}j_{1}}X_{i_{2}j_{1}}X_{i_{2}j_{2}}X_{i_{3}j_{2}}\cdots X_{i_{k}j_{k}}X_{i_{1}j_{k}}),
\end{eqnarray*}
where $k=k_{p}$ satisfies, as $p\to\infty$
\begin{eqnarray*}
\left\{
\everymath{\displaystyle }
\begin{array}{l}
 k/\log{p}\to\infty,\\
 \delta^{1/3}k/\log{p}\to0,\\
 \frac{\delta^{2}\sqrt[4]{p}}{k^{3}}\geq1,
\end{array}
\right.
\end{eqnarray*}
and the summation is taken with respect to $j_{1},j_{2},\ldots,j_{k}$
running over all integers in $\{1,2,\ldots,n\}$ and
$i_{1},i_{2},\ldots,i_{k}$ running over all integers in $\{1,2,\ldots,p\}$
subject to the condition that $i_{1}\neq i_{2},i_{2}\neq
i_{3},\ldots,i_{k}\neq i_{1}$.

In order to get an up bound for $|\mathop{\sum}
EX_{i_1j_1}X_{i_2j_1}\cdots X_{i_kj_k}X_{i_1j_k}|$, we need to construct a graph for given
$i_1,\ldots,i_k$ and $j_1,\ldots,j_k$, as in~\cite{Ger,b2}
and~\cite{b3}. We follow the presentation in~\cite{b3} and
\cite{b2} to introduce some fundamental concepts associated with the
graph.

For the sequence $(i_1,i_2,\ldots,i_k)$ from $\{1,2,\ldots,p\}$ and the
sequence $(j_1,j_2,\ldots,j_k)$ from $\{1,2,\ldots,n\}$, we define a
directed graph as follows. Plot two parallel real lines, referred to
as \emph{I-line} and \emph{J-line}, respectively. Draw
$\{i_1,i_2,\ldots,i_k\}$ on the \emph{I-line}, called \emph{I-vertices}
and draw $\{j_1,j_2,\ldots,j_k\}$ on the \emph{J-line}, known as
\emph{J-vertices}. The vertices of the graph consist of the
\emph{I-vertices} and \emph{J-vertices}. The edges of the graph are
$\{e_1,e_2,\ldots,e_{2k}\}$, where for $a=1,\ldots,k$,
$e_{2a-1}=i_aj_a$ are called the column edges and
$e_{2a}=j_ai_{a+1}$ are called row edges with the convention that
$i_{2k+1}=i_1$. For each column edge $e_{2a-1}$, the vertices
$i_a$ and $j_a$ are called the ends of the edge $i_aj_a$ and
moreover $i_a$ and $j_a$ are, respectively, the initial and the
terminal of the edge $i_aj_a$. Each row edge $e_{2a}$
starts from the vertex $j_b$ and ends with the vertex $i_{b+1}$.

Two vertices are said to coincide if they are both in the
\emph{I-line} or both in the \emph{J-line} and they are identical.
That is $i_a=i_b$ or $j_a=j_b$. Readers are also reminded that the
vertices $i_a$ and $j_b$ are not coincident even if they have the
same value because they are in different lines. We say that two
edges are coincident if two edges have the same set of
ends.\looseness=-1

\begin{figure}

\includegraphics{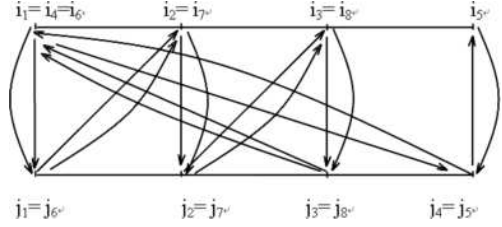}

\caption{An example of \emph{W-graph}.}
\label{fig1}\end{figure}

The graph constructed above is said to be a \emph{W-graph} if each
edge in the graph coincides with at least one other edge. See Figure~\ref{fig1}
for an example of a \emph{W-graph}.


Two graphs are said to be isomorphic if one becomes another by an
appropriate permutation on $\{1,2,\ldots,p\}$ of \emph{I-vertices} and
an appropriate permutation on $\{1,2,\ldots,n\}$ of \emph{J-vertices}.
A \emph{W-graph} is called a canonical graph if $i_a\leq\max\{
i_1,i_2,\ldots,i_{a-1}\}+1 $ and $ j_a\leq\max\{j_1,j_2,\ldots,j_{a-1}\}+1$
with $i_1=j_1=1$,
where $a=1,2,\ldots,k$.


In the canonical graph, if $i_{a+1}=\max\{i_1,i_2,\ldots,i_a\}+1$, then
the edge $j_ai_{a+1}$ is called a row innovation and if $j_a=\max\{
j_1,j_2,\ldots,j_{a-1}\}+1$,
then the edge $i_aj_a$ is called a column innovation.
Apparently, a row innovation and a column innovation, respectively,
lead to a new I-vertex and a new J-vertex except the first column
innovation $i_1j_1$ leading to a new I-vertex $i_1$ and a new J-vertex $j_1$.

We now classify all edges into three types, $T_1$, $T_3$ and $T_4$. Let
$T_1$ denote the set of all innovations including row innovations and
column innovations.
We further distinguish the column innovations as follows. An edge
$i_aj_a$ is called a $T_{11}$ edge if it is a column innovation
and the edge $j_ai_{a+1}$ is a row innovation;
An edge $i_bj_b$ is referred to as a $T_{12}$ edge if it is a
column innovation but $j_bi_{b+1}$ is not a row innovation.
An edge $e_j$ is said to be a $T_3$ edge if there is an innovation edge
$e_i,i<j$ so that $e_j$ is the first one to coincide with
$e_i$. An edge is called a $T_4$ edge if it does not belong to a $T_1$
edge or $T_3$ edge. The first appearance of a $T_4$ edge
is referred to as a $T_2$ edge. There are two kinds of $T_2$ edges: (a)
the first appearance of an edge that coincides with a $T_3$ edge,
denoted by $T_{21}$ edge; (b) the first appearance of an edge that is
not an innovation, denoted by $T_{22}$ edge.

We say that an edge $e_{i}$ is single up to the edge $e_j,j\geq i$, if
it does not coincide with any other edges among $e_1,\ldots,e_j$ except
itself. A $T_3$ edge $e_i$ is said to be regular if there are more than
one innovations with a vertex equal to the initial vertex of $e_{i}$
and single up to $e_{i-1}$, among the edges $\{e_1,\ldots,e_{i-1}\}$.
All other $T_3$ edges are called irregular $T_3$ edges.


Corresponding to the above classification of the edges, we introduce
the following notation and list some useful facts.
\begin{enumerate}[6.]
\item Denote by $l$ the total number of innovations.
\item Let $r$ be the number of the row innovations. Moreover, let $c$
denote the column innovations. We then have $r+c=l$.
\item Define $r_1$ to be the number of the $T_{11}$ edges. Then
$r_1\leq r$ by the definition of a $T_{11}$ edge. Also, the number of
the $T_{12}$ edges is $l-r-r_1$.
\item Let $t$ be the number of the $T_{2}$ edges. Note that the number
of the $T_3$ edges is the same as the number of the innovations and
there are a total of $2k$ edges in the graph. It follows that the
number of
the $T_4$ edges is $2k-2l$. On the other hand, each $T_2$ edge is
also a $T_4$ edge. Therefore, $t\leq2k-2l $.
\item Define $\mu$ to be the number of $T_{21}$ edges. Obviously, $\mu
\leq t$. The number of $T_{22}$ edge is then $t-\mu$. Since each
$T_{21}$ edge coincides with one innovation, we let
$n_{i},i=1,2,\ldots,\mu$, denote the number of $T_{4}$ edges which
coincide with the $i$th such innovation, $n_{i}\geq0$.
\item Let $\mu_1$ be the number of $T_{21}$ edges which do not
coincide with the other $T_4$ edges. That is $\mu_{1}= \#\{
i\dvt n_{i}=1,i=1,2,\ldots,\mu\}$, where $\#\{\cdot\}$ denotes the
cardinality of the set $\{\cdot\}$.
\item Let $m_{j},j=1,2,\ldots,t-\mu$, denote the number of $T_{4}$ edges
which coincide with and include the $j$th $T_{22}$ edge. Note that
$m_{j}\geq2$.\vadjust{\goodbreak}
\end{enumerate}

We now claim that
%
\begin{eqnarray}
\label{a9}
Etr(B_{p}^{k}) & \leq&\bigl(2\sqrt{np}\bigr)^{-k}\sum
E(X_{i_{1}j_{1}}X_{i_{2}j_{1}}\cdots X_{i_{k}j_{k}}X_{i_{1}j_{k}})\nonumber\\
&=& \bigl(2\sqrt{np}\bigr)^{-k}\sum{'}\sum{''}\sum
{'''}\sum_{*}E(X_{i_{1}j_{1}}X_{i_{2}j_{1}}\cdots
X_{i_{k}j_{k}}X_{i_{1}j_{k}})\nonumber
\\[-12pt]
\\[-4pt]
& \leq&
\bigl(2\sqrt{np}\bigr)^{-k}\sum_{l=1}^{k}\sum_{r=1}^{l}\sum_{r_1=0}^{r}\sum
_{t=0}^{2k-2l}\sum_{\mu=0}^{t}\sum_{\mu_1=0}^{\mu}\sum
_{*}\binom{k}{r}\binom{r}{r_1}\binom{k-r_1}{l-r-r_1}\binom
{2k-l}{l}\nonumber\\
&&{}  \times
k^{3t}(t+1)^{6k-6l}\bigl(\delta\sqrt[4]{np}\bigr)^{2k-2l-2t+\mu_{1}}
p^{r+1}n^{l-r},
\nonumber
\end{eqnarray}
where the summation $\sum{'}$ is with respect to different
arrangements of three types of edges at the $2k$ different
positions, the summation $\sum{''}$ over different canonical
graphs with a given arrangement of the three types of edges for $2k$
positions, the third summation $\sum{'''}$ with respect to
all isomorphic graphs for a given canonical graph and the last
notation $\sum_{*}$ denotes the constraint that $i_{1}\neq
i_{2},i_{2}\neq i_{3},\ldots,i_{k}\neq i_{1}$.

\renewcommand\theenumi{(\roman{enumi})}
\renewcommand\labelenumi{\theenumi}
Now, we explain why the above estimate is true:
\begin{enumerate}[(viii)]
\item The factor $(2\sqrt{np})^{-k}$ is obvious.
\item If the graph is not a \emph{W-graph}, which means there is a single
edge in the graph, then the mean of the product of $X_{ij}$
corresponding to this graph is zero (since $EX_{11}=0$). Thus, we have
$l\leq k$. Moreover, the facts that $r\leq l$, $r_1\leq r$, $t\leq
2k-2l$, $\mu\leq t$ and $\mu_1\leq\mu$
are easily obtained from the fact $1$ to the fact $7$ listed before.
\item There are at most $\binoma{k}{r}$ ways to choose $r$ edges out of
the $k$ row edges to be the $r$
row innovations. Subsequently, we consider how to select the column
innovations. Observe that the definition of $T_{11}$ edges,
there are $\binoma{r}{r_1}$ ways to select $r_1$ row innovations out of
the total $r$ row innovations so that the edge before each
such $r_1$ row innovations is a $T_{11}$ edge, column innovation.
Moreover, there are at most $\binoma{k-r_1}{l-r-r_1}$ ways to choose
$l-r-r_1$ edges
out of the remaining $k-r_1$ column edges to be the $l-r-r_1$ $T_{12}$
edges, the remaining column innovations.
\item Given the position of the $l$ innovations, there are at most
$\binoma{2k-l}{l}$ ways to select $l$ edges out of the $2k-l$ edges
to be $T_{3}$ edges. And the rest positions are for the $T_{4}$
edges. Therefore, the first summation $\sum{'}$ is bounded by
$\sum_{l=1}^{k}\sum_{r=1}^{l}\sum
_{r_1=0}^{r}\binoma{k}{r}\binoma{r}{r_1}\binoma{k-r_1}{l-r-r_1}\binoma
{2k-l}{l}$.
\item By definition, each innovation (or each irregular $T_3$ edges) is uniquely
determined by the subgraph prior to the innovation (or the irregular
$T_3$). Moreover, by Lemma 3.2 in~\cite{b2} for each regular $T_3$
edge, there are at most $t+1$ innovations so that the regular $T_3$
edge coincides with one of them and by Lemma 3.3 in~\cite{b2} there
are at most $2t$ regular $T_3$ edges. Therefore, there are at most
$(t+1)^{2t}\leq(t+1)^{2(2k-2l)}$ ways to draw the regular $T_3$
edges.
\item Once the positions of the innovations and the $T_3$ edges are
fixed there are at most $\binoma{(r+1)c}{t}\leq\binoma{k^{2}}{t}\leq
k^{2t}$ ways to arrange the $t$ $T_{2}$ edges, as there are $r+1$ \emph
{I-vertices} and $c$ \emph{J-vertices}. After $t$ positions
of $T_2$ edges are determined there are at most $t^{2k-2l}$ ways to
distribute $2k-2l$ $T_4$ edges among the $t$ positions. So there are
at most $k^{2t}\cdot t^{2k-2l}$ ways to arrange $T_{4}$ edges. It
follows that $\sum{''}$ is bounded by
$\sum_{t=0}^{2k-2l}(t+1)^{2(2k-2l)}k^{2t}\cdot t^{2k-2l}$.

\item The third summation $\sum
{'''}$ is bounded by $n^{c}p^{r+1}$ because the number of
graphs in the isomorphic class for a given graph is
$p(p-1)\cdots(p-r)n(n-1)\cdots(n-c+1)$.

\item Recalling the definitions of $l,r,t,\mu,\mu_{1},n_{i},m_{i}$,
we have
%
\begin{equation}EX_{i_{1}j_{1}}X_{i_{2}j_{1}}\cdots X_{i_{k}j_{k}}X_{i_{1}j_{k}}
=
(EX_{11}^{2})^{l-\mu}\Biggl(\prod_{i=1}^{\mu}EX_{11}^{n_{i}+2}\Biggr)\Biggl(\prod
_{i=1}^{t-\mu}EX_{11}^{m_{i}}\Biggr),\label{a10}
\end{equation}
where $\sum_{i=1}^{\mu}n_{i}+\sum_{i=1}^{t-\mu}m_{i}=2k-2l$. Without
loss of generality, we suppose $n_{1}=n_{2}=\cdots =n_{\mu_{1}}=1$ and
$n_{\mu_{1}+1},\ldots,n_{\mu}\geq2$ for convenience. It is easy to
check that
\begin{eqnarray*}
E|X_{11}^{s}| \leq\left\{
\everymath{\displaystyle }
\begin{array}{l@{ \qquad }l}
M\bigl(\delta\sqrt[4]{np}\bigr)^{s-4}, & \mbox{if } s \geq4,
M=\operatorname{max}\{EX_{11}^{4},|EX_{11}^{3}|\},\\
\bigl(\delta\sqrt[4]{np}\bigr)^{s-2}, & \mbox{if } s\geq2.
\end{array}
\right.
\end{eqnarray*}
Thus, (\ref{a10}) becomes
%
\begin{eqnarray}
\label{a11}
&& |EX_{i_{1}j_{1}}X_{i_{2}j_{1}}\cdots X_{i_{k}j_{k}}X_{i_{1}j_{k}}|\nonumber\\
&& \quad  \leq
\sum_{\mu=0}^{t}\sum_{\mu_1=0}^{\mu}|EX_{11}^{3}|^{\mu
_{1}}|EX_{11}^{4}|^{t-\mu_{1}}\bigl(\delta\sqrt[4]{np}\bigr)^{\sum_{i=\mu
_{1}+1}^{\mu}n_{i}-2(\mu-\mu_{1})}\bigl(\delta\sqrt[4]{np}\bigr)^{\sum
_{i=1}^{t-\mu}m_{i}-2(t-\mu)}\nonumber
\\[-8pt]
\\[-8pt]
&& \quad  \leq\sum_{\mu=0}^{t}\sum_{\mu_1=0}^{\mu}M^{t}\bigl(\delta\sqrt
[4]{np}\bigr)^{2k-2l-2t+\mu_{1}}\nonumber\\
&& \quad  \leq\sum_{\mu=0}^{t}\sum_{\mu_1=0}^{\mu}k^{t}\bigl(\delta\sqrt
[4]{np}\bigr)^{2k-2l-2t+\mu_{1}}, \qquad
\mbox{when $k$ is large enough}.
\nonumber
\end{eqnarray}
\end{enumerate}

The above points regarding the $T_2$ edges are discussed for $t>0$,
but they are still valid when $t=0$ with the convention that $0^0=1$
in the term $t^{2k-2l}$, because in this case there are only $T_1$
edges and $T_3$ edges in the graph and thus $l=k$.

Consider the constraint $\sum_{*}$ now. Note that for each
$T_{12}$ edge, say $i_aj_a$, it is a column innovation, but the next
row edge $j_ai_{a+1}$ is not a row innovation. Since $i_{a+1}\neq
i_a$, the edge $j_{a}i_{a+1}$ cannot coincide with the edge
$i_aj_a$. Moreover, it also doesn't coincide with any edges before
the edge $i_aj_a$ since $j_a$ is a new vertex. So $j_ai_{a+1}$ must
be a $T_{22}$ edge. Thus, the number of the $T_{12}$ edges cannot
exceed the number of the $T_{22}$ edges. This implies $l-r-r_1\leq
t-\mu$. Moreover, note that $\mu_1\leq\mu$. We then have
%
\begin{eqnarray}\label{a12}
&& n^{-k/2}p^{-k/2}n^{l-r}p^{r+1}(np)^{k/2-l/2-t/2+\mu_1/4}\nonumber\\
&& \quad  = (n/p)^{l/2}\cdot n^{-r-t/2+\mu_1/4}p^{r+1-t/2+\mu_1/4}\\
&& \quad  \leq\Biggl(\sqrt{\frac{p}{n}}\Biggr)^{r-r_1}\cdot p^{-t/2}p.\nonumber
\end{eqnarray}
We thus conclude from (\ref{a9}) and (\ref{a12}) that
%
\begin{eqnarray}\label{a13}
Etr(B_{p}^{k}) &\leq&
2^{-k}\sum_{l=1}^{k}\sum_{r=1}^{l}\sum_{r_1=0}^{r}\sum
_{t=0}^{2k-2l}\sum_{\mu=0}^{t}\sum_{\mu_1=0}^{\mu}\binom
{k}{r}\binom{r}{r_1}\binom{k-r_1}{l-r-r_1}\binom{2k-l}{l}\nonumber
\\[-8pt]
\\[-8pt]
&&{}  \times
 \Biggl(\sqrt{\frac{p}{n}}
\Biggr)^{r-r_1}p^{-t/2}pk^{3t}(t+1)^{6k-6l}\delta^{2k-2l-2t+\mu_{1}}.
\nonumber
\end{eqnarray}

Moreover, we claim that
%
\begin{eqnarray}\label{a14}
\hspace*{-15pt}&& p \biggl[2^{-k}\binom{k}{r} \biggr]\Biggl [\binom{r}{r_1} \Biggl(\sqrt
{\frac{p}{n}} \Biggr)^{r-r_1} \Biggr]\biggl [\binom
{k-r_1}{l-r-r_1}\delta^{l-r-r_1} \biggr]  \nonumber\\
\hspace*{-15pt}&& \quad {}  \times\biggl [\binom{2k-l}{l} \biggl(\frac{\sqrt{p}\delta
^3}{k^3} \biggr)^{-t}(t+1)^{6k-6l}\delta^{2k-2l} \biggr]\delta
^{-(l-r-r_1)+3t-(2k-2l)}\cdot\delta^{2k-2l-2t+\mu_1} \\
\hspace*{-15pt}&& \quad  \leq p^2 \Biggl(1+\sqrt{\frac{p}{n}} \Biggr)^k(1+\delta)^k
\biggl(1+\frac{24^3k^3\delta}{\log^3{p}} \biggr)^{2k}.  \nonumber
\end{eqnarray}
Indeed, the above claim is based on the following five facts.
\renewcommand\thelonglist{(\arabic{longlist})}
\renewcommand\labellonglist{\thelonglist}
\begin{longlist}[(5)]
\item$2^{-k}\binoma{k}{r}\leq
2^{-k}\mathop{\sum}_{r=0}^{k}\binoma{k}{r}=1$.
\item $\binoma{r}{r_1} (\sqrt{\frac{p}{n}}
)^{r-r_1}=\binoma{r}{r-r_1} (\sqrt{\frac{p}{n}} )^{r-r_1}
\leq\mathop{\sum}_{s=0}^{r}\binoma{r}{s} (\sqrt{\frac
{p}{n}} )^{s}
=  (1+\sqrt{\frac{p}{n}} )^r
\leq (1+\sqrt{\frac{p}{n}} )^k$.
\item$\binoma{k-r_1}{l-r-r_1}\delta^{l-r-r_1} \leq\mathop{\sum
}_{s=0}^{k-r_1}\binoma{k-r_1}{s}\delta^s
= (1+\delta)^{k-r_1}
\leq(1+\delta)^k$.
\item By the fact that $\binoma{2k-l}{l}\leq\binoma{2k}{2l}$, and the
inequality $a^{-t}(t+1)^b\leq a (\frac{b}{\log{a}} )^b$,
for $a>1$, $b>0$, $t>0$ and $\frac{\delta^2\sqrt{p}}{k^3}\geq\sqrt[4]{p}$,
we have
\begin{eqnarray*}
\binom{2k-l}{l} \biggl(\frac{\sqrt{p}\delta^3}{k^3}
\biggr)^{-t}(t+1)^{6k-6l}\delta^{2k-2l}  &\leq&\binom{2k}{2l}\frac{\sqrt{p}\delta^3}{k^3} \biggl(\frac
{6k-6l}{\log{(\sqrt{p}\delta^3/k^3)}} \biggr)^{6k-6l}\delta
^{2k-2l}\\
&\leq& p\binom{2k}{2l} \biggl(\frac{24k}{\log{p}}
\biggr)^{6k-6l}\delta^{2k-2l}\\
&\leq& p\binom{2k}{2l} \biggl(\frac{24^3k^3\delta}{\log^3{p}}
\biggr)^{2k-2l}\\
&\leq& p\sum_{s=0}^{2k}\binom{2k}{s} \biggl(\frac{24^3k^3\delta}{\log
^3{p}} \biggr)^{2k-s}\\
&=& p \biggl(1+\frac{24^3k^3\delta}{\log^3{p}} \biggr)^{2k}.
\end{eqnarray*}
\item When $p$ is large enough, $\delta^{-(l-r-r_1)+3t-(2k-2l)}\cdot
\delta^{2k-2l-2t+\mu_1}=\delta^{t-(l-r-r_1)}\cdot\delta^{\mu
_1}\leq1$, since $\delta\to0$ and $l-r-r_1\leq t-\mu$.
\end{longlist}

Summarizing (\ref{a13}) and (\ref{a14}), we obtain that
\begin{eqnarray*}
Etr(B_p^k)& \leq&\sum_{l=1}^{k}\sum_{r=1}^{l}\sum_{r_1=0}^{r}\sum
_{t=0}^{2k-2l}\sum_{\mu=0}^{t}\sum_{\mu_1=0}^{\mu} p^2
\Biggl(1+\sqrt{\frac{p}{n}} \Biggr)^k(1+\delta)^k \biggl(1+\frac
{24^3l^3\delta}{\log^3{p}} \biggr)^{2k}\\
& \leq&8k^6p^2 \Biggl(1+\sqrt{\frac{p}{n}} \Biggr)^k(1+\delta)^k
\biggl(1+\frac{24^3l^3\delta}{\log^3{p}} \biggr)^{2k}\\
& = & \Biggl((8k^6)^{1/k}p^{2/k}\Biggl(1+\sqrt{\frac{p}{n}}\Biggr)(1+\delta)
\biggl(1+\frac{24^3k^3\delta}{\log^3{p}} \biggr)^2 \Biggr)^k\\
& \leq&\eta^k,
\end{eqnarray*}
where $\eta$ is a constant satisfying $1<\eta<1+\epsilon$. Here the
last inequality uses the facts below:
\renewcommand\thelonglist{(\roman{longlist})}
\renewcommand\labellonglist{\thelonglist}
\begin{enumerate}[(iii)]
\item$(p^{2})^{1/k}\to1$, because $k/\log{p}\to\infty$,
\item$(8k^{6})^{1/k}\to1$, because $k\to\infty$,
\item$ (1+\sqrt{\frac{p}{n}} )\to1$, because $p/n\to0$,
\item$(1+\delta)\to1$, because $\delta\to0$,
\item$\frac{24^3\cdot k^{3}\delta}{\log^{3}{p}}\to0$, because
$\frac{\delta^{1/3}k}{\log{p}}\to0.$
\end{enumerate}
It follows that
\[
P\bigl(\lambda_{\mathrm{max}}(\bbB_{p})>1+\epsilon\bigr)\leq\biggl(\frac{\eta}{1+\epsilon
}\biggr)^{k}=\mathrm{o}(p^{-\ell})
\]
since $k/\log{p}\to\infty$ and $\frac{\eta}{1+\epsilon}<1$. The proof
is complete.

\section{\texorpdfstring{Proof of Theorem \protect\ref{theo3}}{Proof of Theorem 3}}
Note that
%
\begin{equation}
\bbS_1=\bbS-\bar{\mathbf s}\bar{\mathbf s}'.\label{a15}
\end{equation}
By the Fan inequality~\cite{Fan},
\[
\sup_{x}|F^{\bbA_{p1}}(x)-F^{\bbA_{p}}(x)|\leq\frac{1}{p}.
\]
Thus from theorem in~\cite{b1}, we see that
\[
F^{\bbA_{p1}}(x)\stackrel{\mathrm{a.s.}}\longrightarrow F(x),
\]
specified in the introduction. It follows that
\[
\liminf_{p\rightarrow\infty} \lambda_{\max}(\bbA
_{p1})\geq1.
\]
Let $\bbz$ be a unit vector. In view of (\ref{a15}), we obtain
\[
\bbz'\bbA_{p1}\bbz=\bbz'\bbA_{p}\bbz-\frac{1}{2}\sqrt{\frac
{n}{p}}\bbz'\bar{\mathbf s}\bar{\mathbf s}'\bbz\leq\bbz'\bbA_{p}\bbz,
\]
which implies that
\[
\lambda_{\max}(\bbA_{p1})\leq\lambda_{\max}(\bbA_{p}).
\]
This, together with Theorem~\ref{theo1}, finishes the proof of
Theorem~\ref{theo3}.

\section{\texorpdfstring{Proof of Theorem \protect\ref{theo4}}{Proof of Theorem 4}}

Theorem~\ref{theo4} follows from Theorem~\ref{theo3} and the fact
that
\[
\|\bbS_2-\Sigma\|=\|\Sigma^{1/2}(\bbS_1-\bbI_p)\Sigma^{1/2}\|\leq
\|\bbS_1-\bbI_p\|\|\Sigma\|.
\]

\section*{Acknowledgement}
G.M. Pan was partially supported by a Grant M58110052 at the Nanyang Technological University, and by a grant \# ARC 14/11 from Ministry of Education, Singapore.

%

\printhistory

\end{document}